\documentclass[11pt]{amsart}
\usepackage{amssymb,amsthm,amsmath}
\usepackage{color}

\newcommand{\bt}{\begin{Theorem}}
\newcommand{\et}{\end{Theorem}}
\newcommand{\bi}{\begin{itemize}}
\newcommand{\ei}{\end{itemize}}
\newcommand{\bea}{\begin{eqnarray}}
\newcommand{\ba}{\begin{array}}
\newcommand{\eea}{\end{eqnarray}}
\newcommand{\ea}{\end{array}}

\newenvironment{pf}[1][]{%
	\vskip 3mm
	\noindent
	\ifthenelse{\equal{#1}{}}%
	{{\slshape Proof. }}%
	{{\slshape #1.} }%
}%
{\qed\bigskip}

\usepackage{hyperref}
\usepackage{ifthen}
\usepackage{graphicx}
\nonstopmode
 \numberwithin{equation}{section}

\newtheorem{Definition}{Definition}[section]
\newtheorem{thm}[Definition]{Theorem}
\newtheorem{Lemma}[Definition]{Lemma}

\newtheorem{cor}{Corollary}[section]

\newtheorem*{theoremA*}{Theorem A}
\newtheorem*{theoremB*}{Theorem B}
\newtheorem*{Proofofmainthm*}{Proof of main theorem}
\newcommand{\be}{\begin{equation}}
\newcommand{\ee}{\end{equation}}
\newcommand{\bvb}{\begin{verbatim}}

\allowdisplaybreaks[2]

\newcommand{\evb}{\end{verbatim}}

\newcommand{\x}{\mathfrak X}%

\allowdisplaybreaks[1]

\begin{document}
\baselineskip16pt

\author{Shyam Swarup Mondal} 

\address{Shyam Swarup Mondal  \endgraf Department of Mathematics
	\endgraf IIT Guwahati
	\endgraf Guwahati, Assam, India.} 
\email{mondalshyam055@gmail.com}



\date{\today}

\title[Pseudo differential operators on  Homogeneous Tree]
{Boundedness and Nuclearity of pseudo-differential operators on Homogeneous Trees}

\begin{abstract}
	 Let $\x$ be the  homogeneous tree  of degree $q+1~(q\geq2)$.  In this article, we present symbolic criteria  for  boundedness    of pseudo differential operators on  homogeneous tree $\x$. A sufficient condition for weak type $(p, q)$ boundedness is also given.  We present  necessary and sufficient conditions on the symbols $\sigma$ such that the corresponding pseudo-differential operators $T_\sigma$ from $L^{p_1} (\x)$ into $L^{p_2} (\x)$ to be nuclear for $1\leq p_1, p_2<\infty$.  We show that the adjoint  of the nuclear operator from  $L^{p_2'}(\x)$ into  $L^{p_1'}(\x)$ is again a nuclear operator.   Explicit formulas for the symbol of the adjoint and product of the nuclear pseudo differential operators  on $\x$ are given.
\end{abstract}

\maketitle
\tableofcontents 
\section{Introduction}
The theory of pseudo-differential operators is currently a  very active brand  in  modern mathematics,  mostly for its connections with mathematical physics and  others disciplines like  theory of partial differential equations, harmonic analysis,  black-hole physics, quantum field theory, geometry, time-frequency analysis, micro-local analysis,  imaging, and computations \cite{Hor}. 
Kahn and Nirenberg  \cite{Niren} first introduced pseudo-differential operators   in the mid-1960s  and later  used by H\"ormander \cite{Hor}  to study the problems of geometry and of partial differential equations. Pseudo-differential operators   are extensively  studied by several authors \cite{  Cardona, CardonaKumar1,  CardonaKumar, DK,  Aparajita, Ghaemi, Vish, ssm, Ruz2, RuzT, Wong2} on different classes of groups in different contexts. 

In particular, the definition of pseudo-differential operators on $ \mathbb{ Z}$ follows from the standard idea of starting from the representation of functions
using the inverse of the Fourier transform.  A suitable theory of pseudo-differential operators without symbolic calculus on  $\mathbb{ Z}^n$  has been studided by several authors can be seen in  \cite{Shala3, Delgado, bot}. More recently, a global symbolic calculus on $\mathbb{ Z}^n$ has been introduced and constructed by  Botchway, Kabiti and Ruzhansky \cite{bot}.  Shahla Molahajloo has studied $L^p$-boundedness and $L^p$-compactness of pseudo-differential operators on  $\mathbb{ Z }$ in  \cite{Shala3} (also see \cite{  Carlos}). Boundedness, compactness and the Hilbert-Schmidt property for the $\mathbb{ Z}$-operators  studied by Catana in \cite{cat2}. A more generalization,  $L^p$-boundedness and $L^p$-compactness of multilinear pseudo-differential operators on $\mathbb{ Z}^n$ and $\mathbb{T}^n$ is given in  \cite{cat3, CardonaKumar}.


One of the essential  attraction of weak $L^p$ space is that the subject is sufficiently concrete and yet the spaces have a fine structure of importance for applications in interpolation theory, harmonic analysis, probability theory, functional analysis, and   the Hardy-Littlewood maximal function 
\cite{Bennett}. Weak $L^p$ spaces are function spaces which are closely related to $L^p$ spaces.    
Rodriguez in \cite{Carlos} proved  that $(p, q)$-weak operators on $\mathbb{ Z }^n$ with $p > q \geq  1$ are bounded on $L^p(\mathbb{ Z }^n).$  Weak type multilinear operators on $\mathbb{Z}^n$ has been studied by Catana \cite{cat3}. Later a Weak type (1,1) bounds for a class of operators with discrete kernel has been studied by  Cardona \cite{Cardona1}. The author  also  characterise the weak $(1, p)$-inequality, $1 < p < \infty$, for multipliers on $\mathbb{ Z }^n$ in \cite{cardona22}.

In general, due to Lidskii \cite{lid}, for a trace class operators in Hilbert spaces, the operator trace is given by the sum of the eigenvalues of the operator counted with multiplicities. This property is called the Lidskii formula. However,  in general, the Lidskii formula does not hold in Banach spaces. Grothendieck first proved that for  $2/3$-nuclear operators, the operator trace in Banach spaces agrees with the sum of all the eigenvalues with multiplicities counted \cite{Groth1, Groth2} and this known as the Grothendieck–Lidskii formula. Therefore, for this reason, the concept of nuclear operators are important to us.  The initiative of finding sufficient condition nuclearity and the ${2}/{3}$-nuclearity of pseudo-differential operators on the circle $\mathbb{S}^1$ and on the lattice $\mathbb{ Z}$ has been started by Delgado and Wong \cite{Delgado}. 
Characterizations for nuclear operators on $\mathbb{Z}$ in terms of the decomposition of symbol through Fourier transform has been studied by Jamalpour Birgani in \cite{jamal1}. Nuclearity of pseudo-differential operators for arbitrary compact and Hausdorff groups can be seen in  \cite{Ghaemi}. Recently, the author in \cite{Shyam} discussed nuclearity of $\mathbb{ Z }$-operators related to finite measure spaces.  Nuclearity of pseudo-differential operators on  different class of groups were broadly studied by several authors \cite{Ghasemi,  RuzT, DR, DR18, DRTk} in different context.

A homogeneous tree $\x$ of degree $q+1~(q\geq2)$  is a connected graph with no loops, in which every vertex is adjacent to $q+1$ other vertices (see \cite{FTC,FTP2}). The tree $\x$ carries a natural measure, a counting measure, and a natural distance $d,$ viz, $d(x,y)$ is the number of edges between vertices $x$ and $y$.
Masson  in \cite{ Massan} studied and introduced a general  symbol classes and associated pseudo differential operators on homogeneous tree $\x$. The author  proved that these  operator are   $L^2(\x)$ to  $L^2(\x)$  bounded and gave explicit formulas for adjoint and product. 

In this paper our main goal is to  study  boundedness  and Hilbert-Schmidt property  of pseudo-differential operator on the homogeneous tree $\x.$  We also present  symbolic criteria of $r$-Schatten class operators defined on $L^2(\x)$.  The relation between weak continuity and continuity of a linear operator defined on $L^p(\x)$ is also given. We present a characterization of  nuclear pseudo-differential operator on $\x$ with application to adjoints and products. We give necessary and sufficient conditions on the symbols such that the corresponding pseudo-differential operator on $\x$ to be a nuclear operator. Further, we provide necessary and sufficient conditions on the symbols to guarantee that the adjoints and products of pseudo-differential operators are nuclear. 

The paper is structured in five small sections apart from the introduction. Section \ref{sec2} provides some basic Fourier analysis results  on the homogeneous tree from \cite{Massan}. In section \ref{sec3},  we obtain a necessary and sufficient condition on symbols for which these operators on $\x$ are in the class of Hilbert-Schmidt operators.
We  give a sufficient condition for $L^2$-boundedness  of  pseudo-differential operators on $\x$.  A relation between weak continuity and continuity of a linear operator defined on $L^p(\x)$ given in section \ref{sec4}.   We begin section \ref{sec5}  by presnting a  characterization of nuclear pseudo-differential operator on  $L^p(\x)$.  We also find symbols for the adjoint and product of nuclear pseudo-differential operators on $\x$ and give a characterization for self-adjointness and normality.


\section{Harmonic analysis on Homogeneous Trees}\label{sec2}
We begin this section by recapitulating some basic and important results of harmonic analysis on the homogeneous tree. For more details  about  homogeneous tree see  \cite{CS, Massan, CMS, FTC, FTP2}.

A {\it homogeneous tree} $\x$ of degree $q+1$ is a connected graph with no loops, in which every vertex is adjacent to $q+1$ other vertices.
We assume $q\geq2$.  We fix an arbitrary reference point $o$ in $\x$ and  $|x|$ will denote the distance from $o$ of the vertex $ x$. The tree $\x$ carries a natural measure, a counting measure, and a natural distance $d,$ viz, $d(x,y)$ is the number of edges between the vertices $x$ and $y$. We denote $[x, y]$,  the unique geodesic path joining two arbitrary vertices $x$ and $y$.

An infinite geodesic ray $\omega$ in $\x$ is a one-sided sequence $\{\omega_{n}: n=0,1,2\ldots \}$ where the $\omega_n$ are
in $\x$.    The boundary $\Omega$ of the tree $\x$ is the set of all infinite geodesic rays starting at $o$. It is also the set of equivalence classes of  one-sided  geodesics, where two half geodesics $\{\omega_{n}: n=0,1,2\ldots  \}$ and $\{y_{n}: n=0,1,2\ldots \}$ are equivalent if they meet at infinity, that is if there exists $k, N \in \mathbb{N}$ such that for all $n \geq N, \omega_{n}=y_{n+k}.$ The one-sided  geodesic starting at a point $x$ and equivalent to $\omega$ will be denoted by $[x, \omega)$.

The height functions $h_{\omega}: \mathfrak{X} \rightarrow \mathbb{Z}$ defined for all $\omega=\{x_{0}, x_{1}, \ldots \}$ in $\Omega$ with $x_{0}=o$ by
$$
h_{\omega}(x)=\lim _{m \rightarrow+\infty}\left(m-d\left(x, x_{m}\right)\right).
$$
These functions are normalized so that $h_{\omega}(o)=0 .$ Another definition of the height functions can be give in this way, for every $x \in \mathfrak{X}$ and $\omega \in \Omega,$  we denote $c(x, \omega)$ be the confluence point of $x$ and $\omega $ i.e.,   the last point lying on $[o, \omega)$ in the geodesic path $[o, x]$ then
$$
h_{\omega}(x)=2|c(x, \omega)|-|x|.
$$
For $x\in \x,$ let $E_i(x)=\{\omega\in \Omega : |c(x,\omega)|=i\}$ for all $|x|\geq j\geq0$.  
We use  the  notation $E(x)$ for the set $E_{|x|}(x)$. Then there exist a unique $K$-invariant probability measure $\nu$ such that $$\nu(E_j(x))=\dfrac{q}{(q+1)q^j},\text{  where }j\leq|x|.$$ 
Therefore $(\Omega,\mathcal{M},\nu)$ is a measure space with  $\sigma$ algebra $\mathcal M$ being generated by the set $\{E(x): x\in \x\}$. $\mathcal M_{n}$ be the sub-algebra of $\mathcal M$ generated by the sets $E(x),~|x|\leq n.$ The conditional expectation or  $ \mathcal{E}_{n}(F)$ of an integrable function $F$ on $(\Omega, \mathcal{M}, \nu)$ relative to the sub-algebra $\mathcal{M}_{n}$ is defined by 
$$ \mathcal{E}_{n}(F)(\omega)=\frac{1}{\nu(E_{n}(\omega))}\int\limits_{E_{n}(\omega)}F(\omega')d\nu(\omega').$$
This  $\mathcal{E}_{n}$  are called the family of averaging operators.
\subsection{Change of reference point}
One can define all the previous substance  with respect to another reference point.  For a new reference point \(x_{0}\) we would write \(c^{x_{0}}(x, \omega), E_{i}^{x_{0}}(x) \mathcal{E}_{n}^{x_{0}} \ldots\) to denote the change of reference point. Nonetheless, we will write \(h_{\omega}(x)-h_{\omega}\left(x_{0}\right)\) in place of  \(h_{\omega}^{x_{0}}(x)\). One  also note the Radon-Nikodym derivative $$
\frac{\mathrm{d} \nu_{y}}{\mathrm{d} \nu_{x}}(\omega)=q^{h_{\omega}(y)-h_{\omega}(x)}
$$  for all $x, y \in \mathfrak{X}$ and $\omega \in \Omega $.
\subsection{The Fourier-Helgason transform}
Let $f: \mathfrak{X} \rightarrow \mathbb{C}$  be a function of finite support.  The Fourier-Helgason transform
of $f$  is denoted by $\mathcal{H} f$ defined by  
$$
\mathcal{H} f(\omega, s)=\sum_{y \in \x} f(y) q^{(\frac{1}{2}+i s) h_{\omega}(y)}
$$
for every $\omega \in \Omega$ and $s \in [0, \tau]$ with $ \tau=\frac{\pi}{\log q}.$  Here \(\log (t)\) is the natural logarithm of \(t \in \mathbb{R} .\)  We  also have the   inversion formula given by 
$$
f(x)=\int_{\Omega} \int_{0}^\tau q^{(\frac{1}{2}-i s) h_{\omega}(x)} \mathcal{H}f(\omega, s) \mathrm{d} \nu(\omega) \mathrm{d} \mu(s)
$$
for every $x\in \x$.
The measure $ \mathrm{d} \mu(s)$ is called the Plancherel measure and is  defined by
$$
\mathrm{d} \mu(s)=c_{P}|c(s)|^{-2},
$$
where \(c_{P}=\frac{q \log (q)}{4 \pi(q+1)}\) and for every \(z \in \mathbb{C} \backslash \tau \mathbb{Z}\) we have 
\(c(z)=\frac{q^{1 / 2}}{q+1} \frac{q^{1 / 2+i z}-q^{-1 / 2-i z}}{q^{i z}-q^{-i z}}\).

Moreover  the following version of Plancherel formula holds:
if \(f\) and \(g\) are two finitely supported functions on \(\mathfrak{X},\) then
$$
\sum_{x \in \mathfrak{X}} f(x) \overline{g(x)}=\int_{\Omega}\int_{0}^\tau \mathcal{H}f (\omega, s) \overline{\mathcal{H}g(\omega, s)} \mathrm{d} \nu(\omega) \mathrm{d} \mu(s).
$$

Note that the Fourier-Helgason transform \(\mathcal{H}\) can be extended to an isometry from \(L^{2}(\x)\) to its image in \(L^{2}(\Omega \times \mathbb{T})\).
Now we are in a position  to define pseudo-differential operator on homogeneous tree  using the  inversion formula.
\begin{Definition}
	Let \(\sigma: \mathfrak{X} \times \Omega \times[0, \tau] \rightarrow \mathbb{C}\) be a measurable function. Then for every finitely supported function \(f\) on $\x$, we define  \(T_{\sigma} f\) formally  by 
	\begin{align*}
	T_{\sigma} f(x)&=\int_{\Omega} \int_{0}^\tau  q^{(\frac{1}{2}-i s) h_{\omega}(x)} \sigma(x, \omega, s)\mathcal{H}{f}(\omega, s) \mathrm{d} \nu(\omega) \mathrm{d} \mu(s),~ x\in \x.
	\end{align*}
	We call \(T_{\sigma}\) the pseudo-differential operator  on \(\x\) corresponding to the symbol \(\sigma\).
\end{Definition}
 For more details results related to pseudo differential operator on homogeneous tree and their symbolic calculus we refer \cite{Massan}.

\section{Boundedness }\label{sec3}

This section is devoted to study  boundedness of   pseudo differential operators  on $\x$. We also study  $r$-Schatten-von Neumann  class of pseudo-differential operators on the Hilbert space $L^2(\x)$.

If \(\mathcal{H}\) is a complex Hilbert space, a linear compact operator \(A : \mathcal{H} \rightarrow \mathcal{H}\) belongs to the $r$-Schatten-von Neumann class \(S_{r}(\mathcal{H})\) if
$$
\sum_{n=1}^{\infty}\left(s_{n}(A)\right)^{r}<\infty,
$$where \(s_{n}(A)\) denote the singular values of \(A,\) i.e. the eigenvalues of \(|A|=\sqrt{A^{*} A}\)
with multiplicities counted.

For \(1 \leq r<\infty\),  the class \(S_{r}(\mathcal{H})\) is  a Banach space
endowed with the norm
$$
\|A\|_{S_{r}}=\left(\sum_{n=1}^{\infty}\left(s_{n}(A)\right)^{r}\right)^{\frac{1}{r}}.
$$

For  \(0<r<1\), the $\|\cdot\|_{S_r}$ as above only defines a quasi-norm with respect to which
\(S_{r}(\mathcal{H})\) is complete. An operator belongs to the class \(S_{1}(\mathcal{H})\) is known as {\it Trace class} operator. Also, an operator belongs to   \(S_{2}(\mathcal{H})\) is known as  {\it Hilbert-Schmidt} operator.

We begin this section by  obtaining a necessary and sufficient condition on symbols for which these operators are in the class of Hilbert-Schmidt operators.  Indeed we have the following theorem.

\begin{thm}\label{HS}
	Let $\sigma: \mathfrak{X} \times \Omega \times[0, \tau] \rightarrow \mathbb{C}$ be a measurable function and let \(b: \mathfrak{X} \times \Omega \times[0, \tau]\rightarrow \mathbb{C}\)
	be the measurable function defined by
	$$
	b( x, \omega, s)=\overline{\sigma(x, \omega, s)} q^{( \frac{1}{2} +i s) h_{\omega}(x)},\quad   (x, \omega, s) \in \x \times \Omega \times[0, \tau].
	$$
	Then the  pseudo differential operator \(T_{\sigma}: L^{2}(\x) \rightarrow L^{2}(\x)\) is a Hilbert-Schmidt operator if and only
	if \(b \in L^{2}(\x \times \Omega \times[0, \tau] )\). Moreover 
	$$
	\left\|T_{\sigma}\right\|_{s_{2}}=\|b\|_{L^{2}(\x \times \Omega \times[0, \tau])}.
	$$
\end{thm}
\begin{pf}
	For $x_0\in \x$,  define the function  $f_{x_0}$ on $\x$ by $$  f_{x_0}(y)=\left\{\begin{array}{ll}{0} & {\text{if } x_0\neq y,} \\ {1} & {\text{if } x_0= y.} \end{array}\right. 
	$$ Then $\{f_{x_0}\}_{x_0\in \x}$ be the standard orthonormal basis for $L^2(\x).$ For $x\in \x$ we have 
	\begin{align*}
	(T_\sigma f_{x_0})(x)&=\int_{\Omega}  \int_{0}^\tau q^{( \frac{1}{2} -i s) h_{\omega}(x)} \sigma(x, \omega, s)~ \mathcal{H}f_{x_0}(\omega, s) \mathrm{d} \nu(\omega) \mathrm{d} \mu(s)\\
	&=\sum_{y \in \x} \int_{\Omega}  \int_{0}^\tau q^{( \frac{1}{2} -i s) h_{\omega}(x)} q^{(\frac{1}{2}+i s) h_{\omega}(y)} \sigma(x, \omega, s)~ f_{x_0}(y) \mathrm{d} \nu(\omega) \mathrm{d} \mu(s)\\
	&=\int_{\Omega}  \int_{0}^\tau q^{( \frac{1}{2} -i s) h_{\omega}(x)} q^{(\frac{1}{2}+i s) h_{\omega}(x_0)} \sigma(x, \omega, s) \mathrm{d} \nu(\omega) \mathrm{d} \mu(s)\\
	&= \overline{ \int_{\Omega}  \int_{0}^\tau  q^{(\frac{1}{2}-i s) h_{\omega}(x_0)} \left[ q^{( \frac{1}{2} +i s) h_{\omega}(x)}  \overline{\sigma(x, \omega, s)}\right]  \mathrm{d} \nu(\omega) \mathrm{d} \mu(s)}\\
	&=\overline{ \int_{\Omega}  \int_{0}^\tau  q^{(\frac{1}{2}-i s) h_{\omega}(x_0)} b(x, \omega, s) \mathrm{d} \nu(\omega) \mathrm{d} \mu(s)}\\
	&=\overline{\left( \mathcal{H}^{-1} b(x, \cdot, \cdot)  \right) (x_0)},
	\end{align*}
	where $b(x, \omega, s) =q^{( \frac{1}{2} +i s) h_{\omega}(x)}  \overline{\sigma(x, \omega, s)}$.
	Using  Plancherel formula we get
	\begin{align*}\left\|T_{\sigma}\right\|_{S_{2}}^{2} &=\sum_{x_0 \in \x} \left\|T_{\sigma} f_{x_0}\right\|_{L^{2}(\x)}^{2}\\&=\sum_{x_0 \in \x}  \sum_{x \in \x}  \left|\left(T_{\sigma} f_{x_0}\right)(x)\right|^{2} \\
	&=\sum_{x_0 \in \x}  \sum_{x \in \x}  \left|{\left( \mathcal{H}^{-1} b(x, \cdot, \cdot)  \right) (x_0)}\right|^{2}\\
	&=\sum_{x \in \x}  \sum_{x_0 \in \x}  \left|{\left( \mathcal{H}^{-1} b(x, \cdot, \cdot)  \right) (x_0)}\right|^{2}\\
	&=\sum_{x \in \x}    \left\|{ \mathcal{H}^{-1} b(x, \cdot, \cdot)}  \right\|_{L^{2}(\x)}^{2}\\
	&=\sum_{x \in \x}  \int_{\Omega}  \int_{0}^\tau \left | b(x, \omega, s)  \right |^{2} \mathrm{d} \nu(\omega) \mathrm{d} \mu(s)=\|b\|_{L^{2}(\x \times \Omega \times[0, \tau])}^2.
	\end{align*}
\end{pf}

	The following lemma is a consequence of the definition of Schatten classes (see \cite{Ruz}) which is needed to obtain our  results, $r$-Schatten-von Neumann  class of pseudo-differential operators on the Hilbert space $L^2(\x).$
\begin{Lemma}\label{100}  Let \(A : \mathcal{H} \rightarrow \mathcal{H}\) be a linear compact operator. Let \(0<r, t<\infty .\) Then
	\(A \in S_{r}(\mathcal{H})\) if and only if \(|A|^{\frac{r}{t}} \in S_{t}(\mathcal{H}) .\) Moreover, \(\|A\|_{S_{r}}^{r} =\||A|^{\frac{r}{t}} \|_{S_{t}}^{t}\).
\end{Lemma}
Using the above lemma, we have a characterization of a  pseudo-differential operator on $L^2(\x)$ to be a Schatten class operator. The proof follows from Lemma \ref{100} with $t=2$ and Theorem \ref{HS}.
\begin{cor}
	Let $T : L^2 (\x) \rightarrow L^2 (\x )$ be a continuous linear operator with  symbol \(\sigma\) on $\mathfrak{X} \times \Omega \times[0, \tau].$  Then $T\in S_{r}\left(L^{2}(G/H)\right)  $ if and only if 
	$$ \sum_{x \in \x}  \int_{\Omega}  \int_{0}^\tau \left | \sigma_1(x, \omega, s)  \right |^{2} \mathrm{d} \nu(\omega) \mathrm{d} \mu(s)   <\infty,$$ where $\sigma_1(x, \omega, s)$ is the symbol of the operator $|T|^{\frac{r}{2}}$.
\end{cor}

The next theorem gives a simple sufficient condition on the symbols \(\sigma\) such that the  corresponding  pseudo differential operators \(T_{\sigma}: L^{2}(\x) \rightarrow L^{2}(\x)\) is bounded. A sufficient condition in terms of symbol class  is already given in \cite{Massan}. 
\begin{thm}
	Let \(\sigma:\mathfrak{X} \times \Omega \times[0, \tau] \rightarrow \mathbb{C}\) be a measurable function such that  there
	exist a function \(v \in L^{2}(\x)\) and a positive constant \(C\) for which
	$$
	|q^{ \frac{1}{2} h_{\omega}(x)} \sigma(x, \omega, s)| \leq C|v(x)|, \quad x\in \x,
	$$
	and almost all \((\omega, s)\) in \(\Omega \times[0, \tau]\). Then \(T_{\sigma}: L^{2}(\x) \rightarrow L^{2}(\x)\) is a bounded linear operator.
	Moreover,
	$$
	\left\|T_{\sigma}\right\|_{B(L^{2}(\x))} \leq C'\|v\|_{L^{2}(\x)}
	$$ 
	for some constant $C'.$
\end{thm}

\begin{pf}
	Let $f\in L^2(\x).$ Then  using Cauchy-Schwarz inequality and Plancherel formula 
	\begin{align*}\left\|T_{\sigma}f\right\|_{L^{2}(\x)}^2 &=\sum_{x \in \x} \left |T_{\sigma} f(x)\right |^{2}\\
	&=\sum_{x \in \x} \left |\int_{\Omega}  \int_{0}^\tau q^{( \frac{1}{2} -i s) h_{\omega}(x)} \sigma(x, \omega, s)~ \mathcal{H}f(\omega, s) \mathrm{d} \nu(\omega) \mathrm{d} \mu(s)\right |^{2}\\
	&C_1\leq \sum_{x \in \x} \int_{\Omega}  \int_{0}^\tau \left| q^{( \frac{1}{2} -i s) h_{\omega}(x)} \sigma(x, \omega, s)\right|^2   |\mathcal{H}f(\omega, s) |^2\mathrm{d} \nu(\omega) \mathrm{d} \mu(s)\\
	&C_1\leq \sum_{x \in \x} \int_{\Omega}  \int_{0}^\tau \left| q^{ \frac{1}{2} h_{\omega}(x)} \sigma(x, \omega, s)\right|^2   |\mathcal{H}f(\omega, s) |^2\mathrm{d} \nu(\omega) \mathrm{d} \mu(s)\\
	&\leq C_1C^2\sum_{x \in \x} \left|v(x)\right|^2  \int_{\Omega}  \int_{0}^\tau  |\mathcal{H}f(\omega, s) |^2\mathrm{d} \nu(\omega) \mathrm{d} \mu(s)\\
	&= C'^2 \|v\|_{L^2(\x)}^2\|f\|_{L^2(\x)}^2
	\end{align*}
	and the proof is complete.
\end{pf}

Since a Hilbert-Schmidt operator is a compact operator,  using Theorem \ref{HS} we have the following result.
\begin{cor}
	Let $\sigma: \mathfrak{X} \times \Omega \times[0, \tau] \rightarrow \mathbb{C}$ be a measurable function and let \(b: \mathfrak{X} \times \Omega \times[0, \tau]\rightarrow \mathbb{C}\)
	be the measurable function defined by
	$$
	b( x, \omega, s)=\overline{\sigma(x, \omega, s)} q^{( \frac{1}{2} +i s) h_{\omega}(x)},\quad   (x, \omega, s) \in \x \times \Omega \times[0, \tau].
	$$ Let  \(b \in L^{2}(\x \times \Omega \times[0, \tau] )\),  then the  pseudo differential operator \(T_{\sigma}: L^{2}(\x) \rightarrow L^{2}(\x)\) is a compact operator.
\end{cor}

\section{Weak type $L^p$-estimate }\label{sec4}
In this section, we study the relation between weak continuity and continuity of a linear
operator defined on $L^p(\x)$. In this context, we begin with the definition of $(p, q)$-weak continuity for general measure spaces.

\begin{Definition}
	Let $ (X, \mu)$ and $(Y,\nu)$ be two measure spaces and let $T$ be a linear operator
	from $L^p(X)$ into the space of measurable functions from $Y$ to $\mathbb{C}$.  Then $T$ is said to be a weak type $(p, q)$ or     $(p, q)$-weak operator, $1 \leq  q < \infty$, if there exists $C > 0$ such that 
	\begin{align}\label{1001}
	\nu ( \{y\in Y : |T f (y)| > \lambda \}) \leq  \left( C \frac{\|f\|_{L^p(X)}}{\lambda }\right)^q,~\quad \text{for all $\lambda > 0$}
	\end{align}
	i.e. $T$ maps $L^p(X, \mu) $ boundedly into $L^{q, \infty}(Y, \nu)$.
\end{Definition}
We say that \(T\) is of strong type \((p, q)\) if it is bounded from \(L^{p}(X, \mu)\) to \(L^{q}(Y, \rho) .\) 
Let us also remind that if \((X, \mu)\) is a measure space and \(f: X \rightarrow C\) is a measurable function on \(X,\)
then
$$
\|f\|_{L^{p}(X)}^{p}=p \int_{0}^{\infty} \lambda^{p-1} m(\lambda) d \lambda
$$
for all \(0< p<\infty,\) where \(m(\lambda)=\mu(\{x \in X :|f(x)|>\lambda\})\) is the distribution function of \(f\). From onwards \(\vartheta\) will stands for the counting measure defined on $\x.$

\begin{thm}\label{1002}
	Let T be a linear operator from $L^p(\x)$ into the space of functions from \(\x\) to \(\mathbb{C}\). If \(T\) is a \((p, q)\)-weak operator, \(p>q \geq 1\), then
	$$
	T: L^{p}\left(\x \right) \rightarrow L^{p}\left(\x  \right)
	$$
	is bounded. Moreover     $$
	\|T \|_{B(L^{p}\left(\x \right))} \leq 2^{\frac{p-q}{p}} C \left( \frac{p}{p-q}\right)^{\frac{1}{p}}.
	$$
\end{thm}
\begin{pf}
	Let \(p>q \geq 1\).  Let us assume that \(T\) is  a \((p, q)\)-weak operator. Then from (\ref{1001}),
	$$
	\vartheta ( \{x \in \x : |T f (x)| > \lambda \}) \leq  \left( C \frac{\|f\|_{L^p(\x )}}{\lambda}\right)^q
	$$
	for all \(\lambda>0\) and some \(C>0 .\) Let \(f \in L^{p}\left(\x  \right)\) such that  \(\|f\|_{L^{p}\left(\x \right)}=1 .\) Using the definition 
	\begin{align*}
	\|T f\|_{L^{p}\left(\x\right)}^{p}&=p \int_{0}^{\infty} \vartheta ( \{ x \in \x :|(T f)(x)|>\lambda \}) \lambda^{p-1} d \lambda\\
	&\leq p C^{q} \int_{0}^{2 C} \lambda^{p-q-1} d \lambda\\
	&=p C^{q} \frac{(2 C)^{p-q}}{p-q}.
	\end{align*}
	Therefore, $$ \|T f\|_{L^{p}\left(\x\right)}\leq  2^{\frac{p-q}{p}} C \left( \frac{p}{p-q}\right)^{\frac{1}{p}}.$$    
	Thus for  an arbitrary $f\in L^{p}(\x)$ with \(\|f\|_{L^{p}\left(\x \right)} \neq 0\), we get
	$$
	\|T f\|_{L^{p}\left(\x \right)} \leq 2^{\frac{p-q}{p}} C \left( \frac{p}{p-q}\right)^{\frac{1}{p}} \|f\|_{L^{p}\left(\x \right)}.
	$$
\end{pf}

An immediate consequence of Theorem \ref{1002} gives boundedness of    pseudo differential operator from $L^{p}( \x )$ into $L^{p}( \x),  2< p<\infty$. Indeed, we have the following result.
\begin{cor}
	Let \(2<p<\infty\) such that $\frac{1}{p}+\frac{1}{q}=1$. Let  \(\sigma\) be  a measurable function on $\mathfrak{X} \times \Omega \times[0, \tau]$ satisfying
	$$
	\left(\int_{\Omega} \int_{0}^\tau  q^{(\frac{1}{2}-i s) h_{\omega}(x)}q^{(\frac{1}{2}+i s) h_{\omega}(y)}   \sigma(x, \omega, s) \mathrm{d} \nu(\omega) \mathrm{d} \mu(s) \right)_{(x, y) \in\x\times \x }  \in L^{q}\left(\x\times \x\right).
	$$
	Then the corresponding pseudo differential operator \(T_{\sigma}: L^{p}\left(\x\right) \rightarrow L^{p}\left(\x\right)\) is bounded operator. 
\end{cor}
\begin{pf}
	Let \(2<p<\infty\) such that $\frac{1}{p}+\frac{1}{q}=1$. Let $f\in L^1(\x)$.  Using by H\"older's  inequality we get
	\begin{align*}
	\left| \left(T_{\sigma} f\right)(x) \right|^q&=\left|\int_{\Omega} \int_{0}^\tau  q^{(\frac{1}{2}-i s) h_{\omega}(x)} f(x, \omega, s)\mathcal{H}f(\omega, s) \mathrm{d} \nu(\omega) \mathrm{d} \mu(s)\right|^q\\
	&=\left|\int_{\Omega}  \int_{0}^\tau q^{( \frac{1}{2} -i s)h_{\omega}(x)} \sigma(x, \omega, s)  \left( \sum_{y \in \x} q^{(\frac{1}{2}+i s) h_{\omega}(y)}f(y) \right) ~ \mathrm{d} \nu(\omega) \mathrm{d} \mu(s)\right|^q\\
	&\leq  \left( \sum_{y \in \x} \left| \int_{\Omega}  \int_{0}^\tau q^{( \frac{1}{2} -i s)h_{\omega}(x)} q^{(\frac{1}{2}+i s) h_{\omega}(y)} \sigma(x, \omega, s) \mathrm{d} \nu(\omega) \mathrm{d} \mu(s) \right| |f(y)| \right)^q\\
	&\leq  \left(  \sum_{y \in \x} \left|  \int_{\Omega}  \int_{0}^\tau q^{( \frac{1}{2} -i s)h_{\omega}(x)} q^{(\frac{1}{2}+i s) h_{\omega}(y)} \sigma(x, \omega, s) \mathrm{d} \nu(\omega) \mathrm{d} \mu(s)\right|^q \right) \|f\|_{L^p(\x)}^q
	\end{align*}
	for all $x\in \x.$ Thus
	\begin{align*}
	\|  T_{\sigma} f\|_{L^q(\x)}^q
	&\leq  \left(  \sum_{x \in \x}  \sum_{y \in \x}  \left|  \int_{\Omega}  \int_{0}^\tau q^{( \frac{1}{2} -i s)h_{\omega}(x)} q^{(\frac{1}{2}+i s) h_{\omega}(y)} \sigma(x, \omega, s) \mathrm{d} \nu(\omega) \mathrm{d} \mu(s) \right|^q \right) \|f\|_{L^p(\x)}^q.
	\end{align*}
	Hence \(T_{\sigma}: L^{p}\left(\x \right) \rightarrow L^{q}\left(\x\right)\) is a bounded operator. Moreover 
	\begin{align*}\|  T_{\sigma} \|_{B(L^p(\x), L^q(\x))} &\leq  \left(   \sum_{x \in \x}  \sum_{y \in \x}   \left|   \int_{\Omega}  \int_{0}^\tau q^{( \frac{1}{2} -i s)h_{\omega}(x)} q^{(\frac{1}{2}+i s) h_{\omega}(y)} \sigma(x, \omega, s) \mathrm{d} \nu(\omega) \mathrm{d} \mu(s) \right|^q \right)^{\frac{1}{q}}.
	\end{align*}
	Since $T_{\sigma}$ is a bounded operator it follows that $T_{\sigma}$ is a $(p, q)$-weak operator. By Theorem \ref{1002} and a standard density argument  \(T_{\sigma}: L^{p}\left(\x\right) \rightarrow L^{p}\left(\x\right)\)  bounded for $2<p<\infty$.
\end{pf}

The following theorem is a special case of Theorem \ref{1002}; however, to emphasize the definition of weak continuity, we provide  proof for this theorem.
\begin{thm}
	Let T be a linear operator from $L^p(\x)$ into the space of functions from \(\x\) to \(\mathbb{C}\). If \(T\) is a \((p, q)\)-weak operator for \(q \geq 1\), then
	$$
	T: L^{p}\left(\x \right) \rightarrow L^{r}\left(\x \right)
	$$
	is a bounded operator for all $r>q+1$ . Moreover     
	$$
	\|T \|_{B\left(L^{p}\left(\x\right), L^{r}\left(\x\right)\right ) } \leq  C^{\frac{q}{r}} \left( \sum_{k=2}^\infty \frac{k^q}{(k-1)^r} +(2C)^r \right)^{\frac{1}{r}}.
	$$
\end{thm}
\begin{pf}
	Let \(f \in L^{p}\left(\x \right)\) such that  $\|f\|_{L^{p}\left(\x \right)}=1$ and let $ r>q+1$.  Since  $T$ is a   \((p, q)\)-weak operator, so
	$$
	\vartheta (\left\{x \in \x: |(T f)(x)|>\lambda \right\}) \leq \left( \frac{C}{\lambda}\right) ^{q}
	$$
	for some  $C>0$ and for all $\lambda >0.$ 
	Observe that $$\vartheta (\left\{x \in \x:  \left|\left(Tf\right)(x)\right| > 2 C\right\}) \leq \frac{1}{2^q} <1,$$
	so $\vartheta (\left\{x \in \x:  \left|\left(Tf\right)(x)\right| > 2 C\right\}) =0,$ which implies that $\left|T\left(f\right)(x)\right| \leq 2 C$ for all $x\in \x.$
Thus 
	\begin{align*}
	\x=&\left\{x \in \x: \left|\left(Tf\right)(x)\right| \leq 1\right\} \cup \left\{x \in \x:\left| \left( Tf\right) (x)\right| > 1\right\} \\
	=& \bigcup_{k=2}^{\infty}\left\{x \in \x:\frac{1}{k}<\left|\left(Tf\right)(x)\right| \leq \frac{1}{k-1}\right\}  \cup\left\{x \in \x:  1<\left|\left(Tf\right)(x)\right| \leq 2 C\right\} \end{align*}
	for all \(f \in L^{p}\left( \x \right)\) such that  $\|f\|_{L^{p}\left( \x \right)}=1.$
	Let us denote
	$$A_k=\left\{x \in \x:\frac{1}{k}<\left|\left(Tf\right)(x)\right| \leq \frac{1}{k-1}\right\}, k\geq 2$$ 
	and 
	$$B=\left\{x \in \x:  1<\left|\left(Tf\right)(x)\right| \leq 2 C\right\}.$$
	Then 
	\begin{align*}
	\sum_{ x\in \x} |(Tf)(x)|^r &= \sum_{k=2}^\infty \sum_{ x\in \cup A_k}|(Tf)(x)|^r+ \sum_{ x\in B} |(Tf)(x)|^r\\
	&\leq \sum_{k=2}^\infty \frac{1}{(k-1)^r}\vartheta(A_k)+(2C)^r\vartheta(B)\\
	&\leq  C^q\left(\sum_{k=2}^\infty \frac{k^{q}}{(k-1)^r} +(2C)^r \right )
	\end{align*}
	Since $r>q+1$,  the series  $\sum_{k=2}^\infty \frac{k^{q}}{(k-1)^r}$  is converent. Thus for $f\in L^p(\x)$ we have 
	$$\|Tf\|_{L^r(\x)} \leq  C^{\frac{q}{r}} \left( \sum_{k= 2}^\infty \frac{k^{q}}{(k-1)^r} +(2C)^r \right)^{\frac{1}{r}} \|f\|_{L^p(\x)}.$$
	
\end{pf}

\section{Characterization of nuclear operators}\label{sec5}

We first recall the basic notions of nuclear operators on Banach spaces.
Let    \(T\) be a bounded linear operator from a complex Banach space \(X\) into
another complex Banach space  \(Y\) such that there exist sequences \(\left\{x_{n}^{\prime}\right\}_{n=1}^{\infty}\) in
the dual space \(X^{\prime}\) of \(X\) and \(\left\{y_{n}\right\}_{n=1}^{\infty}\) in \(Y\) such that
$$
\sum_{n=1}^{\infty} \left\|x_{n}^{\prime}\right\|_{X^{\prime}} \left\|y_{n}  \right\|_{Y} <\infty
$$
and $$
T x=\sum_{n=1}^{\infty} x_{n}^{\prime}(x) y_{n}, \quad x \in X.
$$
Then we call \(T : X \rightarrow Y\) a nuclear operator and if \(X=Y,\) then its nuclear trace tr \((T)\)
is given by
$$
\operatorname{tr}(T)=\sum_{n=1}^{\infty} x_{n}^{\prime}\left(y_{n}\right).
$$
It can be proved that the definition of a nuclear operator and the definition of
the trace of a nuclear operator are independent of the choices of the sequences
\(\left\{x_{n}^{\prime}\right\}_{n=1}^{\infty}\) and \(\left\{y_{n}\right\}_{n=1}^{\infty} .\) 

The following theorem proved by Delgado \cite{D2} present a characterization of nuclear operators in term of the conditions on the kernel of the operator.
\begin{thm}\label{11}
	Let \(\left(X_{1}, \mu_{1}\right)\) and \(\left(X_{2}, \mu_{2}\right)\) be  two \(\sigma\)-finite measure spaces. Then \(a\)
	bounded linear operator \(T : L^{p_{1}}\left(X_{1}, \mu_{1}\right) \rightarrow L^{p_{2}}\left(X_{2}, \mu_{2}\right), 1 \leq p_{1}, p_{2}<\infty,\) is
	nuclear if and only if there exist sequences \(\left\{g_{n}\right\}_{n=1}^{\infty}\) in \(L^{p_{1}}\left(X_{1}, \mu_{1}\right)\) and \(\left\{h_{n}\right\}_{n=1}^{\infty}\)
	in \(L^{p_{2}}\left(X_{2}, \mu_{2}\right)\) such that for all \(f \in L^{p_{1}}\left(X_{1}, \mu_{1}\right)\)
	
	$$
	(T f)(x)=\int_{X_{1}} K(x, y) f(y) d \mu_{1}(y), \quad x \in X_{2},
	$$
	where
	$$
	K(x, y)=\sum_{n=1}^{\infty} h_{n}(x) g_{n}(y), \quad x \in X_{2}, y \in X_{1}
	$$ and
	$$
	\sum_{n=1}^{\infty}\left\|g_{n}\right\|_{L^{{p_{1}}^{\prime}}\left(X_{1}, \mu_{1}\right)} \left\|h_{n} \right\|_{L^{p_{2}} \left(X_{2}, \mu_{2}\right)} <\infty.
	$$
\end{thm}
The function \(K\) on \(X_{2} \times X_{1}\) in Theorem \ref{11} is called the kernel of the
nuclear operator \(T : L^{p_{1}}\left(X_{1}, \mu_{1}\right) \rightarrow L^{p_{2}}\left(X_{2}, \mu_{2}\right) .\)

Let \((X, \mu)\) be a \(\sigma\)-finite measure space. Let \(T : L^{p}(X, \mu) \rightarrow L^{p}(X, \mu)\)
\(1 \leq p<\infty,\) be a nuclear operator. Then by Theorem \ref{11}, we can find sequences
\(\left\{g_{n}\right\}_{n=1}^{\infty}\) in \(L^{p^{\prime}}(X, \mu)\) and \(\left\{h_{n}\right\}_{n=1}^{\infty}\) in \(L^{p}(X, \mu)\) such that $$
\sum_{n=1}^{\infty}\left\|g_{n}\right\|_{L^{p^{\prime}}(X, \mu)} \left\|h_{n}\right\|_{L^{p}(X, \mu)} <\infty
$$ 
and for all \(f \in L^{p}(X, \mu)\)
$$
\begin{aligned}(T f)(x) &=\int_{X} K(x, y) f(y) d \mu(y), \quad x \in X, \end{aligned}
$$
where
$$
K(x, y)=\sum_{n=1}^{\infty} h_{n}(x) g_{n}(y), \quad x, y \in X
$$ and it satisfies 
$$\int_{X} |K(x, y)|~d\mu(y)\leq \sum_{n=1}^{\infty } \left\|g_{n}\right\|_{L^{p^{\prime}}(X, \mu)}\left\|h_{n}\right\|_{L^{p}(X, \mu)}.$$

The trace  $\operatorname{tr}(T)$  of   $T : L^{p}(X, \mu)  \rightarrow L^{p}(X, \mu)$  is given by   
\begin{align}\label{88}
\operatorname{tr}(T) =\int_{X} K(x, x) d \mu(x). 
\end{align}

In the next theorem we give a sufficient and necessary condition on the symbol $\sigma$ such that the corresponding pseudo-differential operator  from  $L^{p_{1}}( \x)$ into $ L^{p_{2}}(\x)$ is  nuclear for $1 \leq p_{1}, p_{2}<\infty.$ Indeed we have the following theorem.
\begin{thm}\label{133}
	Let $\sigma:  \mathfrak{X} \times \Omega \times[0, \tau] \rightarrow \mathbb{C}$  be a measurable function.  Then the pseudo differential operator  $T_{\sigma} : L^{p_{1}}( \x ) \rightarrow L^{p_{2}}( \x )$ is  nuclear for \( 1 \leq p_{1}, p_{2}<\infty\)  if and only if there exist sequences \(\left\{g_{k}\right\}_{k=1}^{\infty} \in L^{p_{1}^{\prime}}( \x )\)
	and \(\left\{f_{k}\right\}_{k=1}^{\infty} \in L^{p_{2}}( \x )\) such that
	$$
	\sum_{k=-\infty}^{\infty}\left\|g_{k}\right\|_{ L^{p_{1}^{\prime}}( \x )}\left\|f_{k}\right\|_{ L^{p_{2}}( \x)}<\infty
	$$ and for every $(x, w, s)\in \x \times \Omega \times[0, \tau] $ symbol has  the decomposition
	$$
	\sigma(x, \omega, s) = q^{-( \frac{1}{2} -i s) h_{\omega}(x)}  \sum_{k=-\infty}^{\infty} f_k(x)  \overline{ \mathcal{H}  \overline{g_{k}}(\omega, s)}.$$
\end{thm}

\begin{pf}
	Suppose that \(T_{\sigma} : L^{p_{1}}(\x) \rightarrow L^{p_{2}}(\x)\) is nuclear, where \(1 \leq p_{1}, p_{2}<\infty .\) Then by Theorem \ref{11}, there exist sequences \(\left\{g_{k}\right\}_{k=1}^{\infty}\) in \(L^{p_{1}^{\prime}}( \x)\) and \(\left\{f_{k}\right\}_{k=1}^{\infty}\) in \(L^{p_{2}}( \x)\) such that
	
	$$
	\sum_{k=1}^{\infty}\left\|g_{k}\right\|_{ L^{p_{1}^{\prime}}(\x)} \left\|f_{k}\right\|_{ L^{p_{2}}(\x)}<\infty
	$$
	and for all $f \in L^{p_{1}}( \x),$ we have 
	\begin{align} \label{111}\nonumber
	(T_{\sigma} f)(x)&=\int_{\Omega} \int_{0}^\tau  q^{(\frac{1}{2}-i s) h_{\omega}(x)} f(x, \omega, s)\mathcal{H}f(\omega, s) \mathrm{d} \nu(\omega) \mathrm{d} \mu(s)\\\nonumber
	&=\sum_{y \in \x} \int_{\Omega}  \int_{0}^\tau q^{( \frac{1}{2} -i s) h_{\omega}(x)} q^{(\frac{1}{2}+i s) h_{\omega}(y)} \sigma(x, \omega, s)~ f(y) \mathrm{d} \nu(\omega) \mathrm{d} \mu(s)\\\nonumber
	&= \sum_{y \in \x} \left(  \int_{\Omega}  \int_{0}^\tau q^{( \frac{1}{2} -i s) h_{\omega}(x)} q^{(\frac{1}{2}+i s) h_{\omega}(y)} \sigma(x, \omega, s) \mathrm{d} \nu(\omega) \mathrm{d} \mu(s) \right)  f(y) \\
	&=    \sum_{y \in \x} \sum_{k=-\infty}^{\infty} f_k(x) g_{k}(y) f(y).
	\end{align}
	Let $x_0\in \x$ and the function $f$ is given by $$f(y)=f_{x_0}(y)=\left\{\begin{array}{ll}{0} & {\text{if } x_0\neq y,} \\ {1} & {\text{if } x_0= y.} \end{array}\right. 
	$$
	Then from (\ref{111}) we get 
	
	\begin{align}\label{122}\nonumber
	\sum_{k=-\infty}^{\infty} f_k(x) g_{k}(x_0)&=\int_{\Omega}  \int_{0}^\tau q^{( \frac{1}{2} -i s) h_{\omega}(x)} q^{(\frac{1}{2}+i s) h_{\omega}(x_0)} \sigma(x, \omega, s) \mathrm{d} \nu(\omega) \mathrm{d} \mu(s)\\\nonumber
	&= \overline{ \int_{\Omega}  \int_{0}^\tau  q^{(\frac{1}{2}-i s) h_{\omega}(x_0)} \left[ q^{( \frac{1}{2} +i s) h_{\omega}(x)}  \overline{\sigma(x, \omega, s)}\right]  \mathrm{d} \nu(\omega) \mathrm{d} \mu(s)}\\
	&=\overline{ \int_{\Omega}  \int_{0}^\tau  q^{(\frac{1}{2}-i s) h_{\omega}(x_0)} b(x, \omega, s) \mathrm{d} \nu(\omega) \mathrm{d} \mu(s)},
	\end{align}
	where $b(x, \omega, s) =q^{( \frac{1}{2} +i s) h_{\omega}(x)}  \overline{\sigma(x, \omega, s)}$
	and so
	\begin{align}\label{00001}
	\overline{\left( \mathcal{H}^{-1} b(x, \cdot, \cdot)  \right) (x_0)}
	&=  \sum_{k=-\infty}^{\infty} f_k(x) g_{k}(x_0).
	\end{align}
	Now,
	\begin{align*}
	q^{( \frac{1}{2} -i s) h_{\omega}(x)}  \sigma(x, \omega, s)&=\overline{b(x, \omega, s) }\\ &=  \sum_{x_0\in \x}  q^{(\frac{1}{2}-i s) h_{\omega}(x_0)} \overline{\mathcal{H}^{-1} b(x, \cdot, \cdot) (x_0)}\\
	&= \sum_{x_0\in \x} \left(  \sum_{k=-\infty}^{\infty} f_k(x) g_{k}(x_0)\right)   q^{(\frac{1}{2}-i s) h_{\omega}(x_0)}\\
	&=  \sum_{k=-\infty}^{\infty} f_k(x) \left( \sum_{x_0\in \x}  g_{k}(x_0)  q^{(\frac{1}{2}-i s) h_{\omega}(x_0)} \right) \\
	&= \sum_{k=-\infty}^{\infty} f_k(x)  \overline{ \mathcal{H}  \overline{g_{k}} (\omega, s)}.
	\end{align*}
	Therefore,
	\begin{align*}
	\sigma(x, \omega, s) = q^{-( \frac{1}{2} -i s) h_{\omega}(x)}  \sum_{k=-\infty}^{\infty} f_k(x)  \overline{ \mathcal{H}  \overline{g_{k}}(\omega, s)},~~ (x, \omega, s) \in   \x  \times \Omega \times [0, \tau].
	\end{align*}

	Conversely, suppose that there exist sequences \(\left\{g_{k}\right\}_{k=1}^{\infty}\) in \(L^{p_{1}^{\prime}}(\x)\) and \(\left\{f_{k}\right\}_{k=1}^{\infty}\)
	in \(L^{p_{2}}(\x )\) such that
	$$
	\sum_{k=-\infty}^{\infty}\left\|g_{k}\right\|_{ L^{p_{1}^{\prime}}(\x )} \left\|f_{k}\right\|_{ L^{p_{2}}(\x )}<\infty
	$$
	and  for every $(x, w, s)\in \x \times \Omega \times[0, \tau] $
	$$
	\sigma(x, \omega, s) = q^{-( \frac{1}{2} -i s) h_{\omega}(x)}  \sum_{k=-\infty}^{\infty} f_k(x)  \overline{ \mathcal{H}  \overline{g_{k}}(\omega, s)}.
	$$
	Then, for all \(f \in L^{p_{1}}(\x)\)
	\begin{align*}
	(T_{\sigma} f)(x)&=\int_{\Omega}  \int_{0}^\tau q^{( \frac{1}{2} -i s) h_{\omega}(x)} \sigma(x, \omega, s)~ \mathcal{H}f(\omega, s) \mathrm{d} \nu(\omega) \mathrm{d} \mu(s)\\
	&=\int_{\Omega}  \int_{0}^\tau  \left( \sum_{k=-\infty}^{\infty} f_k(x)  \overline{ \mathcal{H}  \overline{g_{k}}(\omega, s)} \right) \mathcal{H}f(\omega, s) \mathrm{d} \nu(\omega) \mathrm{d} \mu(s)\\
	&=\int_{\Omega}  \int_{0}^\tau  \left( \sum_{k=-\infty }^{\infty} f_{k}(x)  \sum_{y \in \x} {g_k(y)} q^{(\frac{1}{2}-i s) h_{\omega}(y)}   \right) \mathcal{H}f(\omega, s) \mathrm{d} \nu(\omega) \mathrm{d} \mu(s)\\
	&= \sum_{y \in \x}  \sum_{k=-\infty}^{\infty}f_{k}(x)  {g_k(y)}  \int_{\Omega}  \int_{0}^\tau q^{(\frac{1}{2}-i s) h_{\omega}(y)}   \mathcal{H}f(\omega, s) \mathrm{d} \nu(\omega) \mathrm{d} \mu(s)\\
	&= \sum_{y \in \x} \sum_{k=-\infty}^{\infty} f_{k}(x)  {g_k(y)}   f(y)
	\end{align*}
	for all \(x \in \x .\) Therefore by Theorem \ref{11},   \( T_{\sigma} : L^{p_{1}}(\x) \rightarrow L^{p_{2}}( \x) \) is  a nuclear operator.
\end{pf}

In order to find the trace of a nuclear  operator from  $L^{p}( \x )$ into $L^{p}( \x )$, in the next theorem, we give another characterization of nuclear  operator from  $ L^{p_{1}}( \x ) $ into  $L^{p_{2}}( \x ), 1 \leq p_{1}, p_{2}<\infty.$

\begin{thm}\label{144}
	Let \( \sigma \) be a measurable  function on \(  \mathfrak{X} \times \Omega \times[0, \tau]\). Then the pseudo differential operator  $T_{\sigma} : L^{p_{1}}( \x ) \rightarrow L^{p_{2}}( \x )$ is  nuclear for \( 1 \leq p_{1}, p_{2}<\infty\)  if and only if there exist sequences \(\left\{g_{k}\right\}_{k=1}^{\infty} \in L^{p_{1}^{\prime}}( \x )\)
	and \(\left\{f_{k}\right\}_{k=1}^{\infty} \in L^{p_{2}}( \x )\) such that
	$$
	\sum_{k=-\infty}^{\infty}\left\|g_{k}\right\|_{ L^{p_{1}^{\prime}}( \x )}\left\|f_{k}\right\|_{ L^{p_{2}}( \x)}<\infty
	$$ and for all $x, y\in \x$ we have
	$$
	\int_{\Omega}  \int_{0}^\tau q^{( \frac{1}{2} -i s) h_{\omega}(x)} q^{( \frac{1}{2} +i s) h_{\omega}(y)} \sigma(x, \omega, s)~ \mathrm{d} \nu(\omega) \mathrm{d} \mu(s)=\sum_{k=-\infty}^{\infty} f_{k}(x)  {g_k(y)}.
	$$
\end{thm}

\begin{pf}
	Let  \(T_{\sigma} : L^{p_{1}}(\x ) \rightarrow L^{p_{2}}(\x )\) is nuclear operator for $1 \leq p_{1}, p_{2}<\infty$. By  Theorem \ref{11} and equation (\ref{122}) there exist sequences \(\left\{g_{k}\right\}_{k=1}^{\infty} \in L^{p_{1}^{\prime}}(\x )\)
	and \(\left\{f_{k}\right\}_{k=1}^{\infty} \in L^{p_{2}}(\x )\) such that
	$$
	\sum_{k=-\infty}^{\infty}\left\|g_{k}\right\|_{ L^{p_{1}^{\prime}}(\x )}\left\|f_{k}\right\|_{ L^{p_{2}}(\x  )}<\infty
	$$ and for every $x, y \in \x$,
	\begin{align*}
	\sum_{k=-\infty}^{\infty} f_k(x) g_{k}(y)=\int_{\Omega}  \int_{0}^\tau q^{( \frac{1}{2} -i s) h_{\omega}(x)} q^{(\frac{1}{2}+i s) h_{\omega}(y)} \sigma(x, \omega, s) \mathrm{d} \nu(\omega) \mathrm{d} \mu(s).
	\end{align*}
	Conversely, suppose there exist sequences  \(\left\{g_{k}\right\}_{k=1}^{\infty} \in L^{p_{1}^{\prime}}(\x )\)
	and \(\left\{f_{k}\right\}_{k=1}^{\infty} \in L^{p_{2}}(\x )\) such that
	$$
	\sum_{k=-\infty}^{\infty}\left\|g_{k}\right\|_{ L^{p_{1}^{\prime}}(\x )}\left\|f_{k}\right\|_{ L^{p_{2}}(\x  )}<\infty
	$$  and for every $x, y \in \x$,
	\begin{align*}
	\sum_{k=-\infty}^{\infty} f_k(x) g_{k}(y)=\int_{\Omega}  \int_{0}^\tau q^{( \frac{1}{2} -i s) h_{\omega}(x)} q^{(\frac{1}{2}+i s) h_{\omega}(y)} \sigma(x, \omega, s) \mathrm{d} \nu(\omega) \mathrm{d} \mu(s).
	\end{align*}
	Now, for all \(f \in L^{p_{1}}(\x)\)
	\begin{align*}
	\left(T_{\sigma} f\right)(x)&=\int_{\Omega}  \int_{0}^\tau q^{( \frac{1}{2} -i s) h_{\omega}(x)} \sigma(x, \omega, s)~ \mathcal{H}f(\omega, s) \mathrm{d} \nu(\omega) \mathrm{d} \mu(s)\\
	&=\sum_{y \in \x} \int_{\Omega}  \int_{0}^\tau q^{( \frac{1}{2} -i s) h_{\omega}(x)} q^{(\frac{1}{2}+i s) h_{\omega}(y)} \sigma(x, \omega, s)~ f(y) \mathrm{d} \nu(\omega) \mathrm{d} \mu(s)\\
	&=\sum_{y \in \x} \sum_{k=-\infty}^{\infty} f_k(x) g_{k}(y) f(y)
	\end{align*}
	for all \(x \in \x .\) Therefore by Theorem \ref{11},   \( T_{\sigma} : L^{p_{1}}(\x) \rightarrow L^{p_{2}}(\x)\) is  a nuclear operator.
\end{pf}

An immediate consequence of Theorem \ref{144} gives the trace of a nuclear   of pseudo differential operator from $L^{p}( \x )$ into $L^{p}( \x ),  1\leq p<\infty$. 
\begin{cor}\label{155}
	Let \(T_{\sigma} : L^{p}( \x ) \rightarrow L^{p}( \x )\) be a  nuclear operator for \(1 \leq p<\infty\).
	Then the nuclear  trace $\operatorname{Tr} \left(T_{\sigma}\right)$ of \(T_{\sigma}\) is given by
	$$
	\operatorname{Tr}\left(T_{\sigma}\right)= \sum_{x \in \x}  \int_{\Omega}  \int_{0}^\tau q^{h_{\omega}(x)} \sigma(x, \omega, s) \mathrm{d} \nu(\omega) \mathrm{d} \mu(s).
	$$
	Moreover, $$
	\left\{ \int_{\Omega}  \int_{0}^\tau q^{h_{\omega}(\cdot)} \sigma(\cdot, \omega, s) \mathrm{d} \nu(\omega) \mathrm{d} \mu(s)\right\}   \in L^1(\x).
	$$
\end{cor}
\begin{pf}
	Using trace formula  (\ref{88}) and Theorem  \ref{144} we have 
	\begin{align*} \operatorname{Tr}\left(T_{\tau}\right) 
	&=\sum_{x \in \x} \sum_{k=-\infty}^{\infty} f_k(x) g_{k}(x)\\
	&=\sum_{x \in \x}  \int_{\Omega}  \int_{0}^\tau q^{( \frac{1}{2} -i s) h_{\omega}(x)} q^{(\frac{1}{2}+i s) h_{\omega}(x)} \sigma(x, \omega, s) \mathrm{d} \nu(\omega) \mathrm{d} \mu(s)\\
	&= \sum_{x \in \x}  \int_{\Omega}  \int_{0}^\tau q^{h_{\omega}(x)} \sigma(x, \omega, s) \mathrm{d} \nu(\omega) \mathrm{d} \mu(s).
	\end{align*}
	
	Using H\"older's inequality  
	\begin{align*}
	\sum_{x \in \x} \left|  \int_{\Omega}  \int_{0}^\tau q^{h_{\omega}(x)} \sigma(x, \omega, s) \mathrm{d} \nu(\omega) \mathrm{d} \mu(s)\right| 
	&=\sum_{x \in \x} \left| \sum_{k=-\infty}^{\infty} f_k(x) g_{k}(x)\right|\\
	&\leq \sum_{k=-\infty}^{\infty}\left\|g_{k}\right\|_{ L^{p'}(\x )}\left\|f_{k}\right\|_{ L^{p}(\x  )}<\infty
	\end{align*} 
	and this implies that $$
	\left\{ \int_{\Omega}  \int_{0}^\tau q^{h_{\omega}(\cdot)} \sigma(\cdot, \omega, s) \mathrm{d} \nu(\omega) \mathrm{d} \mu(s)\right\}   \in L^1(\x).
	$$
\end{pf}

\subsection{Adjoints}
Let $T_{ \sigma } : L^{p_{1}}( \x ) \rightarrow L^{p_{2}} (\x), 1\leq p_{1}, p_{2} < \infty,$  be a nuclear operator. In this sub-section we show that the adjoint operator \(T_{ \sigma }^* : L^{p_{2}^{\prime}} ( \x) \rightarrow L^{p_{1}^{\prime}} (\x) \) of \(T_{ \sigma }\) is  also nuclear. We also  present a necessary and sufficient condition on the symbol \({ \sigma }\) so
that the corresponding muclear operator \(T_{ \sigma }\) from \(L^{2}(  \x  )\) into \(L^{2}(  \x )\) to be self-adjoint.

\begin{thm}    \label{adjoint}
	Let \( \sigma \) be a measurable  function on \( \x \times \Omega \times[0, \tau]\). Let  $T_{\sigma} : L^{p_{1}}( \x ) \rightarrow L^{p_{2}}( \x )$ is  nuclear for \( 1 \leq p_{1}, p_{2}<\infty\). Then $T_{\sigma}^*$,  the adjoint of  $T_{\sigma},$  is also a nuclear operator from  $L^{p_{2}'}(\x)$ into $ L^{p_{1}'}(\x)$ with symbol $\sigma^*$ given by 
	
	$$ {\sigma^*(x, \omega, s)}= q^{-( \frac{1}{2} -i s) h_{\omega}(x)} \sum_{k=-\infty}^{\infty} \overline{\mathcal{H} f_k(\omega, s) }~\overline{g_{k}(x) }, \quad (x, \omega, s)\in \mathfrak{X} \times \Omega \times[0, \tau]$$
	where  \(\left\{g_{k}\right\}_{k=1}^{\infty} \) and \(\left\{f_{k}\right\}_{k=1}^{\infty} \) are two sequences in  \(L^{p_{1}^{\prime}}( \x)\)
	and \(L^{p_{2}}( \x)\) respectively  such that
	$$
	\sum_{k=-\infty}^{\infty}\left\|g_{k}\right\|_{ L^{p_{1}^{\prime}}( \x )}\left\|h_{k}\right\|_{ L^{p_{2}}( \x )}<\infty.
	$$
\end{thm}

\begin{pf}
	For all $f \in L^{p_{1}}( \x)$ and $ g \in L^{p_{2}}( \x)$, we have 
	$$ \sum_{x\in \x } (T_{\sigma}f)(x) \overline{g(x)} =  \sum_{x\in \x  } f(x) \overline{T_{\sigma}^*g(x)}.$$
	Therefore,
	\begin{align} \label{270}\nonumber
	&\sum_{x\in \x} \left (  \sum_{y \in \x}   \int_{\Omega}  \int_{0}^\tau q^{( \frac{1}{2} -i s) h_{\omega}(x)} q^{(\frac{1}{2}+i s) h_{\omega}(y)} \sigma(x, \omega, s) \mathrm{d} \nu(\omega) \mathrm{d} \mu(s) f(y)\right ) \overline{g(x)} \\
	&=  \sum_{x\in \x } f(x) \overline{ \left (  \sum_{y \in \x}   \int_{\Omega}  \int_{0}^\tau q^{( \frac{1}{2} -i s) h_{\omega}(x)} q^{(\frac{1}{2}+i s) h_{\omega}(y)} \sigma^*(x, \omega, s) \mathrm{d} \nu(\omega) \mathrm{d} \mu(s) g(y) \right ) }.
	\end{align}
	\noindent
	For all  \(a, b \in \x\) and for all \(x\in \x\), $f$ and $g$ is defined by 
	$$f(x)=f_{a}(x)=\left\{\begin{array}{ll}{0} & {\text{if } a\neq x,} \\ {1} & {\text{if } a= x,} \end{array}\right. \quad g(x)=g_{b}(x)=\left\{ \begin{array}{ll}{0} & {\text{if } b\neq x,} \\ {1} & {\text{if } b= x.} \end{array}\right. 
	$$
	Therefore from (\ref{270}) we have,
	\begin{align*} 
	&\int_{\Omega}  \int_{0}^\tau q^{(\frac{1}{2}+i s) h_{\omega}(a)} [q^{( \frac{1}{2} -i s) h_{\omega}(b)}  \sigma(b, \omega, s) ]\mathrm{d} \nu(\omega) \mathrm{d} \mu(s) \\
	&= \overline{  \int_{\Omega}  \int_{0}^\tau q^{(\frac{1}{2}+i s) h_{\omega}(b)}[  q^{( \frac{1}{2} -i s) h_{\omega}(a)} \sigma^*(a, \omega, s)] \mathrm{d} \nu(\omega) \mathrm{d} \mu(s)   }
 	\end{align*}
	and thus  for all $a, b\in \x$,
	\begin{align}  \label{280}
	\overline{\left( \mathcal{H}^{-1} \sigma'(b, \cdot, \cdot)  \right) (a)}=
	{\left( \mathcal{H}^{-1} \sigma''(a, \cdot, \cdot)  \right) (b)},
	\end{align} 
	where $\sigma'(b, \omega, s) =q^{( \frac{1}{2} +i s) h_{\omega}(b)}  \overline{\sigma(b, \omega, s)}$ and $\sigma''(a, \omega, s) =q^{( \frac{1}{2} +i s) h_{\omega}(a)}  \overline{\sigma^*(a, \omega, s)}$.
	Since  \(T_{\sigma} : L^{p_{1}}( \x) \rightarrow L^{p_{2}}( \x)\) is  nuclear, from Theorem \ref{133}, there  exist sequences \(\left\{g_{k}\right\}_{k=1}^{\infty} \in L^{p_{1}^{\prime}}( \x )\)
	and \(\left\{f_{k}\right\}_{k=1}^{\infty} \in L^{p_{2}}( \x )\) such that
	$\displaystyle 
	\sum_{k=-\infty}^{\infty}\left\|g_{k}\right\|_{ L^{p_{1}^{\prime}}( \x )}\left\|f_{k}\right\|_{ L^{p_{2}}( \x)}<\infty
	$ and for every $(x, w, s)\in \x \times \Omega \times[0, \tau] $ symbol has  the decomposition
	$$
	\sigma(x, \omega, s) = q^{-( \frac{1}{2} -i s) h_{\omega}(x)}  \sum_{k=-\infty}^{\infty} f_k(x)  \overline{ \mathcal{H}  \overline{g_{k}}(\omega, s)}.$$
	Using the relation (\ref{280}) and the nuclearity of $T_{\sigma}$, we have 
	\begin{align*}
	q^{( \frac{1}{2} +i s) h_{\omega}(a)}  \overline{\sigma^*(a, \omega, s)}&=    \sigma''(a, \omega, s) \\ &=  \sum_{b\in \x}  q^{(\frac{1}{2}+i s) h_{\omega}(b)} (\mathcal{H}^{-1} \sigma''(a, \cdot, \cdot)) (b)\\
	&=\sum_{b\in \x}  q^{(\frac{1}{2}+i s) h_{\omega}(b)}\overline{\left( \mathcal{H}^{-1} \sigma'(b, \cdot, \cdot)  \right) (a)}\\
	&=\sum_{b\in \x}  q^{(\frac{1}{2}+i s) h_{\omega}(b)} \int_{\Omega}  \int_{0}^\tau q^{( \frac{1}{2} +i s_1) h_{\omega_1}(a)} \overline{\sigma'(b, \omega_1, s_1) }\mathrm{d} \nu(\omega_1) \mathrm{d} \mu(s_1) \\
	&=\sum_{b\in \x}  q^{(\frac{1}{2}+i s) h_{\omega}(b)} \int_{\Omega}  \int_{0}^\tau q^{( \frac{1}{2} +i s_1) h_{\omega_1}(a)} q^{( \frac{1}{2} -i s_1) h_{\omega_1}(b)}  \sigma(b, \omega_1, s_1)\mathrm{d} \nu(\omega_1) \mathrm{d} \mu(s_1) \\
	&=\sum_{b\in \x}  q^{(\frac{1}{2}+i s) h_{\omega}(b)} \int_{\Omega}  \int_{0}^\tau q^{( \frac{1}{2} +i s_1) h_{\omega_1}(a)} \sum_{k=-\infty}^{\infty} f_k(b)  \overline{ \mathcal{H}  \overline{g_{k}}(\omega_1, s_1)}\mathrm{d} \nu(\omega_1) \mathrm{d} \mu(s_1) \\
	&=\sum_{b\in \x}  q^{(\frac{1}{2}+i s) h_{\omega}(b)} \sum_{k=-\infty}^{\infty} f_k(b) \int_{\Omega}  \int_{0}^\tau q^{( \frac{1}{2} +i s_1) h_{\omega_1}(a)}   \overline{ \mathcal{H}  \overline{g_{k}}(\omega_1, s_1)}\mathrm{d} \nu(\omega_1) \mathrm{d} \mu(s_1) \\
	&=\sum_{b\in \x}  q^{(\frac{1}{2}+i s) h_{\omega}(b)} \sum_{k=-\infty}^{\infty} f_k(b) g_{k}(a) \\
	&=\sum_{k=-\infty}^{\infty} \sum_{b\in \x}  q^{(\frac{1}{2}+i s) h_{\omega}(b)}  f_k(a) g_{k}(a) \\
	&=\sum_{k=-\infty}^{\infty}  \mathcal{H} f_k(\omega, s) g_{k}(a).
	\end{align*}
	Therefore,  for every $(x, w, s)\in \x \times \Omega \times[0, \tau] $
	$$ {\sigma^*(x, \omega, s)}= q^{-( \frac{1}{2} -i s) h_{\omega}(x)} \sum_{k=-\infty}^{\infty} \overline{\mathcal{H} f_k(\omega, s) }~\overline{g_{k}(x) }.$$
\end{pf}

As a consequence of the above theorem and help of Theorem \ref{133}, we can give a criterion in terms of  symbol such that the corresponding pseudo differential operator  on $L^2$ to  be self-adjoint.

\begin{cor}
	Let \(\sigma\) be a measurable  function on \( \mathfrak{X} \times \Omega \times[0, \tau]\) such that  \(T_{\sigma} : L^{2}( \x) \rightarrow L^{2}( \x)\) is nuclear. Then \(T_{\sigma} : L^{2}( \x) \rightarrow L^{2}( \x)\) is self-adjoint if and only if there exist sequences     \(\left\{g_{k}\right\}_{k=1}^{\infty}\)  and \(\left\{f_{k}\right\}_{k=1}^{\infty}\) in \(L^{2}(\x)\) such that
	
	$$\sum_{k=-\infty }^{ \infty} \left\| h_{k} \right\|_{ L^{2}(\x) } \left \|g_{k} \right \|_{L^{2}(\x)} <\infty,$$
	
	$$    \sum_{k=-\infty}^{\infty} f_k(x)  \overline{ \mathcal{H}  \overline{g_{k}}(\omega, s)}= \sum_{k=-\infty}^{\infty} \overline{\mathcal{H} f_k(\omega, s) }~\overline{g_{k}(x) }$$
	and  $$
	\sigma(x, \omega, s) = q^{-( \frac{1}{2} -i s) h_{\omega}(x)}  \sum_{k=-\infty}^{\infty} f_k(x)  \overline{ \mathcal{H}  \overline{g_{k}}(\omega, s)},$$ for all $ (x, \omega, s)\in \mathfrak{X} \times \Omega \times[0, \tau]    $.
	
\end{cor}
We can give another formula for the adjoints of nuclear pseudo differential operators in terms of symbols. Indeed, we have the following theorem.
\begin{thm}
	Let \(\sigma\) be a measurable  function on \( \mathfrak{X} \times \Omega \times[0, \tau]\) such that  \(T_{\sigma} : L^{2}( \x) \rightarrow L^{2}( \x)\) is nuclear. Then for all $ (y, \omega_1, s_1)\in \mathfrak{X} \times \Omega \times[0, \tau]    $, we have 
	\begin{align*}
 {\sigma^*(y, \omega_1, s_1)}   &=q^{-( \frac{1}{2} -i s_1) h_{\omega_1}(y)}  \sum_{ x\in \x}  q^{( \frac{1}{2} -i s_1) h_{\omega_1}(x)}\\&\times \int_{\Omega}  \int_{0}^\tau   \overline{\sigma(x, \omega, s)} q^{( \frac{1}{2} +i s) h_{\omega}(x)} q^{( \frac{1}{2} -i s) h_{\omega}(y)} \mathrm{d} \nu(\omega) \mathrm{d} \mu(s) .
	\end{align*}
\end{thm}
\begin{pf}
	Let  $T_{\sigma} : L^{p_{1}}( \x ) \rightarrow L^{p_{2}}( \x )$ be a nuclear operator for \( 1 \leq p_{1}, p_{2}<\infty\). Then by Theorem \ref{133} there exist sequences \(\left\{g_{k}\right\}_{k=1}^{\infty} \in L^{p_{1}^{\prime}}( \x )\)
	and \(\left\{f_{k}\right\}_{k=1}^{\infty} \in L^{p_{2}}( \x )\) such that
	$$
	\sum_{k=-\infty}^{\infty}\left\|g_{k}\right\|_{ L^{p_{1}^{\prime}}( \x )}\left\|f_{k}\right\|_{ L^{p_{2}}( \x)}<\infty
	$$ and for every $(x, w, s)\in \x \times \Omega \times[0, \tau] $ symbol has  the decomposition
	$$
	\sigma(x, \omega, s) = q^{-( \frac{1}{2} -i s) h_{\omega}(x)}  \sum_{k=-\infty}^{\infty} f_k(x)  \overline{ \mathcal{H}  \overline{g_{k}}(\omega, s)}.$$
	Let  $y\in \x.$ Then 
	\begin{align*}
	&\int_{\Omega}  \int_{0}^\tau   \overline{\sigma(x, \omega, s)} q^{( \frac{1}{2} +i s) h_{\omega}(x)} q^{( \frac{1}{2} -i s) h_{\omega}(y)} \mathrm{d} \nu(\omega) \mathrm{d} \mu(s) \\&=\int_{\Omega}  \int_{0}^\tau    q^{( \frac{1}{2} -i s) h_{\omega}(y)} \left( \sum_{k=-\infty}^{\infty}\overline{f_k(x) } \mathcal{H}  \overline{g_{k}}(\omega, s)\right)  \mathrm{d} \nu(\omega) \mathrm{d} \mu(s) \\
	&=  \sum_{k=-\infty}^{\infty} \overline{f_k(x) }  \int_{\Omega}  \int_{0}^\tau   q^{( \frac{1}{2} -i s) h_{\omega}(y)} \mathcal{H}  \overline{g_{k}}(\omega, s) \mathrm{d} \nu(\omega) \mathrm{d} \mu(s) \\
	&=   \sum_{k=-\infty}^{\infty} \overline{f_k(x)}~ \overline{g_k(y)}.
	\end{align*}
	Thus, using Theorem \ref{adjoint} we get
	\begin{align*}
	& \sum_{ x\in \x}  q^{( \frac{1}{2} -i s_1) h_{\omega_1}(x)}\int_{\Omega}  \int_{0}^\tau   \overline{\sigma(x, \omega, s)} q^{( \frac{1}{2} +i s) h_{\omega}(x)} q^{( \frac{1}{2} -i s) h_{\omega}(y)} \mathrm{d} \nu(\omega) \mathrm{d} \mu(s) \\&= \sum_{ x\in \x}  q^{( \frac{1}{2} -i s_1) h_{\omega_1}(x)}   \sum_{k=-\infty}^{\infty} \overline{f_k(x)}~ \overline{g_k(y)}\\
	&=  \sum_{k=-\infty}^{\infty}  \left( \sum_{ x\in \x}  q^{( \frac{1}{2} -i s_1) h_{\omega_1}(x)}    \overline{f_k(x)}\right) \overline{g_k(y)}\\
	&= \sum_{k=-\infty}^{\infty}   \overline{(\mathcal{H} f_k)(\omega_1, s_1)}~ \overline{g_k(y)}=   {\sigma^*(y, \omega_1, s_1)}  q^{( \frac{1}{2} -i s_1) h_{\omega_1}(y)} .
	\end{align*}
	Therefore, 	 for every $(y, w_1, s_1)\in \x \times \Omega \times[0, \tau] $ we get 
	\begin{align*}
	{\sigma^*(y, \omega_1, s_1)}   &=q^{-( \frac{1}{2} -i s_1) h_{\omega_1}(y)}  \sum_{ x\in \x}  q^{( \frac{1}{2} -i s_1) h_{\omega_1}(x)}\\&\times \int_{\Omega}  \int_{0}^\tau   \overline{\sigma(x, \omega, s)} q^{( \frac{1}{2} +i s) h_{\omega}(x)} q^{( \frac{1}{2} -i s) h_{\omega}(y)} \mathrm{d} \nu(\omega) \mathrm{d} \mu(s) .
	\end{align*}
\end{pf}

Using the  above theorem, we have  another   criterion in terms of symbol such that the corresponding pseudo differential operator to be self-adjoint.

\begin{cor}
	Let \(\sigma\) be a measurable  function on \( \mathfrak{X} \times \Omega \times[0, \tau]\) such that  \(T_{\sigma} : L^{2}( \x) \rightarrow L^{2}( \x)\) is nuclear. Then \(T_{\sigma} : L^{2}( \x) \rightarrow L^{2}( \x)\) is self-adjoint if and only if 
\begin{align*}
{\sigma(y, \omega_1, s_1)}   &=q^{-( \frac{1}{2} -i s_1) h_{\omega_1}(y)}  \sum_{ x\in \x}  q^{( \frac{1}{2} -i s_1) h_{\omega_1}(x)}\\&\times \int_{\Omega}  \int_{0}^\tau   \overline{\sigma(x, \omega, s)} q^{( \frac{1}{2} +i s) h_{\omega}(x)} q^{( \frac{1}{2} -i s) h_{\omega}(y)} \mathrm{d} \nu(\omega) \mathrm{d} \mu(s) 
\end{align*}
for all $ (y, \omega_1, s_1)\in \mathfrak{X} \times \Omega \times[0, \tau]    $.
\end{cor}

\subsection{Products}
In this subsection, we show that the product of a nuclear pseudo-differential operator with a bounded operator is again a nuclear operator, and as a consequence, we give a necessary and sufficient condition on symbol so that the corresponding nuclear pseudo-differential operator to be normal.

\begin{thm}
	Let $\sigma$ and $ \eta $ be two  measurable function on \(  \x \times \Omega \times[0, \tau]  \) such that  the corresponding  pseudo differential operator $T_{\sigma}:  L^{p}(\x) \rightarrow L^{p}(\x)$ is nuclear and  \(T_{\eta} : L^{p}(\x) \rightarrow L^{p}(\x)\)  is a bounded linear operator for $1\leq p<\infty$. Then $T_{\eta} T_{\sigma}:  L^{p}(\x) \rightarrow L^{p}(\x)$ is a nuclear  operator $T_{\lambda}:  L^{p}(\x) \rightarrow L^{p}(\x)$ with symbol $\lambda $ is given by 
	$$ \lambda(b, \omega, s) = q^{-( \frac{1}{2} -i s) h_{\omega}(b)}   \sum_{x\in \x}     q^{( \frac{1}{2} -i s) h_{\omega}(x)}  {\sigma(x, \omega, s)} {\left( \mathcal{H}^{-1} \eta'(x, \cdot, \cdot)  \right) (b)},$$
	  where $\eta'(x, \omega, s) =q^{( \frac{1}{2} +i s) h_{\omega}(x)}  \overline{\eta^*(x, \omega, s)}, ~(x, \omega, s)\in \x \times \Omega \times[0, \tau] $.
	
\end{thm}
\begin{pf}
	For all $f , g\in  L^{p}(\x)$, we have 
	
	$$\sum_{x\in \x} (T_{\lambda}f)(x) \overline{ {g(x)}}= \sum_{x\in \x} (T_{\eta} T_{\sigma} f)(x) \overline{  g(x)}= \sum_{x\in \x} (T_{\sigma}f)(x) \overline{ T_{\eta^*}  g(x)}.$$
	Therefore,
	\begin{align} \label{250}\nonumber
	&\sum_{x\in \x} \left (  \sum_{y \in \x}   \int_{\Omega}  \int_{0}^\tau q^{( \frac{1}{2} -i s) h_{\omega}(x)} q^{(\frac{1}{2}+i s) h_{\omega}(y)} \lambda(x, \omega, s) \mathrm{d} \nu(\omega) \mathrm{d} \mu(s) f(y)\right ) \overline{g(x)} \\\nonumber
	&=  \sum_{x\in \x } (T_\sigma f)(x) \overline{ \left (  \sum_{y \in \x}   \int_{\Omega}  \int_{0}^\tau q^{( \frac{1}{2} -i s) h_{\omega}(x)} q^{(\frac{1}{2}+i s) h_{\omega}(y)} \eta^*(x, \omega, s) \mathrm{d} \nu(\omega) \mathrm{d} \mu(s) g(y) \right ) }\\\nonumber
	&=  \sum_{x\in \x }    \left (  \sum_{y \in \x}   \int_{\Omega}  \int_{0}^\tau q^{( \frac{1}{2} -i s) h_{\omega}(x)} q^{(\frac{1}{2}+i s) h_{\omega}(y)} \sigma(x, \omega, s) \mathrm{d} \nu(\omega) \mathrm{d} \mu(s) f(y) \right )\\
	&\quad\quad \quad \times \overline{ \left (  \sum_{y \in \x}   \int_{\Omega}  \int_{0}^\tau q^{( \frac{1}{2} -i s) h_{\omega}(x)} q^{(\frac{1}{2}+i s) h_{\omega}(y)} \eta^*(x, \omega, s) \mathrm{d} \nu(\omega) \mathrm{d} \mu(s) g(y) \right ) }
	\end{align}
	\noindent
	For all  \(a, b \in \x\) and for all \(x\in \x\), $f$ and $g$ is defined by 
	$$f(y)=f_{a}(y)=\left\{\begin{array}{ll}{0} & {\text{if } a\neq y,} \\ {1} & {\text{if } a= y.} \end{array}\right. g(y)=g_{b}(y)=\left\{ \begin{array}{ll}{0} & {\text{if } b\neq y,} \\ {1} & {\text{if } b= y.} \end{array}\right. 
	$$
	Therefore from (\ref{250}) we have,
	
	\begin{align*} 
	&\int_{\Omega}  \int_{0}^\tau q^{( \frac{1}{2} -i s) h_{\omega}(b)} q^{(\frac{1}{2}+i s) h_{\omega}(a)} \lambda(b, \omega, s) \mathrm{d} \nu(\omega) \mathrm{d} \mu(s) \\
	&=  \sum_{x\in \x }    \left (   \int_{\Omega}  \int_{0}^\tau q^{( \frac{1}{2} -i s) h_{\omega}(x)} q^{(\frac{1}{2}+i s) h_{\omega}(a)} \sigma(x, \omega, s) \mathrm{d} \nu(\omega) \mathrm{d} \mu(s)  \right )\\
	&\quad\quad \quad \times \overline{ \left ( \int_{\Omega}  \int_{0}^\tau q^{( \frac{1}{2} -i s) h_{\omega}(x)} q^{(\frac{1}{2}+i s) h_{\omega}(b)} \eta^*(x, \omega, s) \mathrm{d} \nu(\omega) \mathrm{d} \mu(s) \right ) }.
	\end{align*}
	Therefore,
	\begin{align*} 
	{\left( \mathcal{H}^{-1} \lambda'(b, \cdot, \cdot)  \right) (a)}=\sum_{x\in \x} {\left( \mathcal{H}^{-1} \sigma'(x, \cdot, \cdot)  \right) (a)}\overline{\left( \mathcal{H}^{-1} \eta'(x, \cdot, \cdot)  \right) (b)}
	\end{align*} 
	where $$\lambda'(b, \omega, s) =q^{( \frac{1}{2} +i s) h_{\omega}(b)}  \overline{\lambda(b, \omega, s)},$$ 
	$$\sigma'(x, \omega, s) =q^{( \frac{1}{2} +i s) h_{\omega}(x)}  \overline{\sigma(x, \omega, s)}$$ and
	$$\eta'(x, \omega, s) =q^{( \frac{1}{2} +i s) h_{\omega}(x)}  \overline{\eta^*(x, \omega, s)}.$$
	So,
	\begin{align*}
	\lambda'(b, \omega, s)
	&=  \sum_{a\in \x}  q^{(\frac{1}{2}+i s) h_{\omega}(a)} (\mathcal{H}^{-1} \lambda'(b, \cdot, \cdot)) (a)\\
	&= \sum_{a\in \x}  q^{(\frac{1}{2}+i s) h_{\omega}(a)} \sum_{x\in \x} \left( \mathcal{H}^{-1} \sigma'(x, \cdot, \cdot)  \right) (a)\overline{\left( \mathcal{H}^{-1} \eta'(x, \cdot, \cdot)  \right) (b)}\\
	&=  \sum_{x\in \x}  \left( \sum_{a\in \x}  q^{(\frac{1}{2}+i s) h_{\omega}(a)} \left( \mathcal{H}^{-1} \sigma'(x, \cdot, \cdot)  \right) (a) \right) \overline{\left( \mathcal{H}^{-1} \eta'(x, \cdot, \cdot)  \right) (b)}\\
	&=  \sum_{x\in \x}   \sigma'(x, \omega, s) \overline{\left( \mathcal{H}^{-1} \eta'(x, \cdot, \cdot)  \right) (b)}.
	\end{align*}
	Therefore for all $(b, \omega, s)\in \x \times \Omega \times[0, \tau] $ we get
	$$ \lambda(b, \omega, s) = q^{-( \frac{1}{2} -i s) h_{\omega}(b)}   \sum_{x\in \x}     q^{( \frac{1}{2} -i s) h_{\omega}(x)}  {\sigma(x, \omega, s)} {\left( \mathcal{H}^{-1} \eta'(x, \cdot, \cdot)  \right) (b)}. $$
	Since $T_{\sigma} $ is nuclear operator and $T_{\eta }$ is bounded operator, by applying Theorem \ref{133}, the proof  will be complete.
\end{pf}

The following corollary gives a criterion in terms of symbol such that the corresponding pseudo differential operator  on $\x$ to be normal.
\begin{cor}
	Let $\sigma$ and $ \eta $ be two  measurable function on \(  \x \times \Omega \times[0, \tau]  \) such that  the corresponding  pseudo differential operators $T_{\tau}:  L^{p}(\x) \rightarrow L^{p}(\x)$ is nuclear and  \(T_{\sigma} : L^{p}(\x) \rightarrow L^{p}(\x)\)  is a bounded linear operator for $1\leq p<\infty$. Then $T_{\sigma} T_{\eta}:  L^{p}(\x) \rightarrow L^{p}(\x)$ is a nuclear  operator $T_{\lambda}:  L^{p}(\x) \rightarrow L^{p}(\x)$ with symbol $\lambda $ is given by 
	$$ \lambda(b, \omega, s) = q^{-( \frac{1}{2} -i s) h_{\omega}(b)}   \sum_{x\in \x}     q^{( \frac{1}{2} -i s) h_{\omega}(x)}  {\eta(x, \omega, s)} {\left( \mathcal{H}^{-1} \sigma'(x, \cdot, \cdot)  \right) (b)},$$
	  where $\sigma'(x, \omega, s) =q^{( \frac{1}{2} +i s) h_{\omega}(x)}  \overline{\sigma^*(x, \omega, s)},  ~(x, \omega, s)\in \x \times \Omega \times[0, \tau] $.
\end{cor}
\begin{cor}
	Let \(\sigma\) be a  function on \( \x \times \Omega \times[0, \tau] \) such that  \(T_{\sigma} : L^{2}(\x) \rightarrow L^{2}(\x)\) is a  nuclear operator. Then \(T_{\sigma} : L^{2}(\x) \rightarrow L^{2}(\x)\) is normal  if and only if  for all $(b, \omega, s)\in \x \times \Omega \times[0, \tau] $
	$$\sum_{x\in \x}     q^{( \frac{1}{2} -i s) h_{\omega}(x)}  {\sigma(x, \omega, s)} {\left( \mathcal{H}^{-1} \sigma'(x, \cdot, \cdot)  \right) (b)} =\sum_{x\in \x}     q^{( \frac{1}{2} -i s) h_{\omega}(x)}  {\sigma(x, \omega, s)} {\left( \mathcal{H}^{-1} \sigma'(x, \cdot, \cdot)  \right) (b)},$$ where $\sigma'(x, \omega, s) =q^{( \frac{1}{2} +i s) h_{\omega}(x)}  \overline{\sigma(x, \omega, s)}$.
\end{cor}


\section{Acknowledgement}

The author thanks   IIT Guwahati for providing financial support. Author also thanks Sumit Kumar Rano  and Vishvesh Kumar for their  fruitful suggestions.

\section{Data availability statement}
The authors confirm that the data supporting the findings of this study are available within the article  and its supplementary materials.

\end{document}